\documentclass[11pt]{amsart}
\usepackage{amssymb}
\usepackage{amsmath}
\usepackage{amsfonts}

\newtheorem{thm}{Theorem}[section]

\newtheorem{defn}[thm]{Definition}

\newtheorem{prob}[thm]{Problem}

\newtheorem{facts}[thm]{Facts}
\newtheorem{remark}[thm]{Remark}

\newtheorem{lemma}[thm]{Lemma}
\newtheorem{proposition}[thm]{Proposition}
\newtheorem{theorem}[thm]{Theorem}
\newtheorem{corollary}[thm]{Corollary}

\theoremstyle{definition}

\newtheorem{example}[thm]{Example}

\begin{document}

\title{Frames, Graphs and Erasures}

\author[B.~G.~Bodmann]{Bernhard G.~Bodmann}

\address{Dept.~of Mathematics, University of Houston, Texas  
77204-3008,
U.S.A.}

\email{bgb@math.uh.edu}

\author[V.~I.~Paulsen]{Vern I.~Paulsen}

\address{Dept.~of Mathematics, University of Houston, Texas  
77204-3008,
U.S.A.}

\email{vern@math.uh.edu}


\begin{abstract}
Two-uniform frames and their use for the coding of vectors
are the main subject of this paper.
These frames are known to be optimal for
handling up to two erasures,
in the sense that they minimize the largest possible error when 
up to two frame coefficients are set to zero. 
Here, we consider various numerical measures for the reconstruction
error associated with a frame when an arbitrary number of the frame coefficients of a vector are 
lost. We derive general error bounds for two-uniform frames when more than 
two erasures occur and apply these to concrete examples. We show that
among the 227 known equivalence classes of two-uniform (36,15)-frames arising from 
Hadamard matrices, there are 5 that give smallest error bounds for up to 8 erasures.
\end{abstract}

\maketitle

%%%%%%%%%%%%%%%%%%%%%%%%%%%%%%%%%%%%%%%%%%%%%%%%%%%%%%%%%%%%%%%%%%%%%%%%%%%%
%%%

% \baselineskip=20pt

%This sets up thm, prop, lem, cor, dfn, rmk the way we want them.

\maketitle
\newcommand{\ms}{\medskip}
\newcommand{\bs}{\bigskip}
\newcommand{\ep}{\epsilon}
\newcommand{\bb}[1]{\mathbb{#1}}
\newcommand{\cl}[1]{\mathcal{#1}}
\newcommand{\what}{\widehat}
\newcommand{\ovl}{\overline}
\newcommand{\ds}{\displaystyle}

\newcommand{\rE}{\tilde{E}}
\def\lam{\lambda}
\def\alg{\text{ alg}}
\def\ball{\text{ ball}}
\def\dist{\text{ dist}}
\def\diag{\text{ diag}}
\def\var{\varphi}
\def\span{\text{ span}}
\def\Im{\text{Im }}
\def\D{\mathbb D}
\def\id{\text{ id}}
\newcommand{\n}[1]{\left\vert#1\right\vert}   %absolute value
\def\bv \bign#1{\bigl\vert#1\bigr\vert} %big absolute value
\newcommand{\V}[1]{\left\Vert#1\right\Vert}   %norm
\def \p#1{\left(#1\right)}         %matched parentheses
\def\ra{\rightarrow}
\def \a#1{\langle#1\rangle}
\def \supnorm#1{\N{#1}_{\infty}}        % || x ||_{\infty}
\def \onenorm#1{\N{#1}_{1}}             % || x ||_1
\def \Xa{\bold X}
\def \R{\mathbb R}
\def\cO{\Cal O}
\def\cX{\Cal X}
\def\cM{\Cal M}
\def\cZ{\Cal Z}
\def\cA{\Cal A}
\def\cN{\Cal N}
\def\cB{\Cal B}
\def\cU{\Cal U}
\def\cH{\Cal H}
\def\cJ{\Cal J}
\def\cI{\Cal I}
\def\cR{\Cal R}
\def\cS{\Cal S}
\def\cK{\Cal K}
\def\be{\beta}
\def\ep{\epsilon}
\def\al{\alpha}
\def\raw{\overset\rightarrow \to w}

\def\ralam{\overset\rightarrow\to  \lambda}
\def\ramu{\overset\rightarrow\to \mu}
\def\raz{\overset\rightarrow \to z}
\def\rax{\overset\rightarrow \to x}
\def\ray{\overset\rightarrow \to y}
\def\raal{\overset\rightarrow\to \alpha}
\def\rabe{\overset\rightarrow\to\beta}
\def\ga{\gamma}
\def\raga{\overset\rightarrow\to\gamma}
\def \Rn{\R^n}                   %fancy R^n
\def \Z{\mathbb Z}
\def\C{\mathbb C}              %fancy Z -- integers
\def\cS{\Cal S}
\def \cC{\Cal C}
\def\cD{\Cal D}
\def \F{\mathbb F}
\def \N{\mathbb N}
\def\cP{\Cal P}
\def\cF{\Cal F}
\def \Sa{\mathbb S}
\def \MIN{\text{MIN}}
\def \MAX{\text{MAX}}

\def\raw{\overset\ra\to w}
\def\FO{\overset\frown\to\otimes}
\def\SO{\underset\smile\to\otimes}
\def\pp{\prime\prime}
\def\Bar{\overline}

%%%%%%%%%%%%%%%%%%%%%%%%%%%%%%%%%%%%%%%%%%%%%%%%%%%%%%%%%%%%%%%%%%%%%%%%%%%%
%%%

\section{Introduction} 

Uniform tight frames are a well-known means for the 
redundant encoding of vectors in terms of their expansion coefficients. 
Such frames were studied in \cite{CK}, 
\cite{GKK} and \cite{GVT} and shown to be optimal in some sense for one 
erasure. In addition, further properties of these frames were developed,
including their robustness to more than one erasure.
In this paper we introduce some measures of how well a frame behaves 
under multiple erasures and 
then seek optimal frames in these contexts.

In an earlier paper \cite{HP}, a family of so-called {\it two-uniform} 
frames was introduced. When they exist, two-uniform
frames were demonstrated to be optimal for one and two erasures. 
Moreover, it was proved that a
frame is two-uniform if and only if it is {\it equiangular} which is a
family of frames that has been studied independently by Thomas
Strohmer and Robert Heath~\cite{St}.
The existence of such frames, over the reals, depends on the existence
of a matrix of $\pm 1's$ which satisfies certain algebraic equations.
Strohmer and Heath made the key discovery that these matrices had been
studied earlier in the graph theory literature and correspond
precisely to the Seidel adjacency matrices of a very special type of
graph.
In this paper, we derive explicit formulas that relate how well two-uniform 
frames behave under erasures to various connectivity problems of the related 
graphs.

This paper is organized as follows. After fixing the notation in Section~2,
we discuss the idea of using frames as codes in Section~3 and introduce a 
a family of numerical measures for the error when the coded information is partially deleted.
Section~4 recalls two-uniform frames as the ones that perform best under
one and two deletions. The construction of such frames is then related to a problem
in graph theory. In Section~5, we derive error bounds from the graph-theoretic formulation.
The general error bounds are illustrated with examples in Section~6.

\section{Preliminaries and Notation}\label{sec:prel}

We begin by recalling the basic definitions and concepts.

\begin{defn}
Let $\cl H$ be a Hilbert space, real or complex, and let $F = \{ 
f_i \}_{i \in 
\mathbb I} \subset \cl H$ be a subset.
We call $F$ a
{\bf frame} for $\cl H$ provided that there are 
two constants $C,D > 0$ such that the norm inequalities
\[ 
 C \cdot \|x\|^2 \leq \sum_{j \in \mathbb I} | \langle x,f_j 
\rangle |^2 \leq 
  D \cdot \|x\|^2 
\] 
hold for every $x \in \cl H$. Here, $\langle\cdot,\cdot\rangle$ denotes the 
inner product of two vectors, by convention conjugate linear in the second entry if $\cl H$
is a complex Hilbert space.

When $C = D = 1$, then we call $F$ a {\bf Parseval frame.} Such frames
are also called {\bf normalized, tight frames,} but Parseval frame is,
perhaps, becoming more standard.

A frame $F$ is called {\bf uniform} or {\bf equal-norm} 
provided there is a constant 
$c$ so that $\|f\| = c$ for all $f \in F.$

The linear map $V : \cl H \to \ell_2(\bb I)$ defined by $$(Vx)_i = 
\langle x,f_i \rangle$$
is called the {\bf analysis operator.} 
When $F$ is a Parseval frame, then $V$ is an isometry 
and its adjoint, $V^*$, 
acts as a left inverse to $V$. 
\end{defn}

For the purposes of this paper we will only be concerned with finite
dimensional Hilbert spaces and frames for these spaces that consist of
finitely many vectors.  When the dimension of $\cl H$ is $k$, then we
will identify $\cl H$ with $\bb R^k$ or $\bb C^k$ depending on whether
we are dealing with the real or complex case, and for notational
purposes regard vectors as columns.

When we wish to refer to either case, then we will denote the 
ground field by $\bb F.$

\begin{defn}
We shall let $\cl F(n,k)$ denote the collection of all Parseval frames
for $\bb F^k$ consisting of $n$ vectors and refer to such a frame as
either a real or complex {\bf (n,k)-frame,} depending on whether or
not the field $\bb F$ is the real numbers or the complex numbers.  Thus, a
uniform $(n,k)$-frame is a uniform Parseval frame for $\bb F^k$ with $n$ vectors.
The ratio $n/k$ we shall refer to as the {\bf redundancy} of the
frame.
\end{defn}

If we identify the analysis operator $V$ of an $(n,k)$-frame 
with an $n \times k$ matrix, using the standard basis, then the 
columns of $V^*$ are the 
frame vectors.

\begin{facts}
  Using some basic operator theory, $F$ is an $(n,k)$-frame if and
  only if the Grammian (or correlation) matrix $VV^* = (\langle
  f_j,f_i \rangle)$ of the frame vectors is a self-adjoint $n \times
  n$ projection of rank $k$.  Moreover, the rank of a projection is
  equal to its trace, so $tr(VV^*)=k.$  Thus, when $F$ is a uniform
  $(n,k)$-frame, each of the diagonal entries
  of $VV^*$ must be equal to $k/n,$ and hence each frame vector $f_j$ must be of length
  $||f_j||=\sqrt{k/n}.$

  Conversely, given an $n \times n$ self-adjoint projection $P$ of
  rank $k$, we can always factor it as $P = VV^*$ with an $n \times k$
  matrix $V$, by choosing an orthonormal basis for the range of $P$ as
  the column vectors of $V$.  It follows that $V^*V = I_k$ and hence
  $V$ is the matrix of an isometry and so corresponds to an
  $(n,k)$-frame. Moreover, if $P= WW^*$ is another factorization of
  $P$, then necessarily $W^*W=I_k$ and there exists a unitary $U$ such 
  that $W^*= UV^*$ and hence
  the two corresponding frames differ by multiplication by this
  unitary.  Thus, $P$ determines a unique unitary equivalence class of
  frames. A self-adjoint projection $P$ corresponds to a uniform 
  $(n,k)$-frame if and only if all of its diagonal entries are $k/n.$
\end{facts}

\begin{defn}
  In the following, we wish to identify certain frames as being
  equivalent.  Given frames $F =\{ f_1, \ldots , f_n \}$ and $G = \{
  g_1, \ldots , g_n \}$, we say that they are {\bf type I equivalent}
  if there exists a unitary (orthogonal, in the real case) matrix $U$
  such that $g_i = Uf_i$ for all $i$. If $V$ and $W$ are the analysis
  operators for $F$ and $G$, respectively, then it is clear that $F$
  and $G$ are type I equivalent if and only if $V=WU$ or equivalently,
  if and only if $VV^* =WW^*.$ Thus, there is a one-to-one
  correspondence between $n \times n$ rank $k$ projections and type I
  equivalence classes of $(n,k)$-frames.

  We say that two frames are {\bf type II equivalent} if they are
  simply a permutation of the same vectors and {\bf type III
    equivalent} if the vectors differ by multiplication with $\pm 1$ in the
  real case and multiplication by complex numbers of modulus one in
  the complex case.

  Finally, we say that two frames are {\bf equivalent} if they belong
  to the same equivalence class in the equivalence relation generated
  by these three equivalence relations.  It is not hard to see that if
  $F$ and $G$ are frames with analysis operators $V$ and $W$,
  respectively, then they are equivalent if and only if $UVV^*U^* =
  WW^*$ for some $n \times n$ unitary $U$ that is the product of a
  permutation and a diagonal unitary (or diagonal orthogonal matrix, 
  in the real case).
\end{defn}

We caution the reader that the equivalence relation that we have just
defined is different from the equivalence relation that is often used,
but it is the one studied in \cite{HP}. In other papers, often frames
$\{ f_i \}$ and $\{ g_i \}$ are called equivalent provided that there
is an invertible operator $T$ such that $Tf_i = g_i$ for all $i$,
which is clearly a much coarser equivalence relation than the one
used here.

\section{Frames and Erasures}

The idea behind treating frames as codes is that, given an original
vector $x$ in $\bb F^k$ and an $(n,k)$-frame with analysis operator
$V$, one regards the vector $Vx\in \bb F^n$ as an encoded version of $x$, 
which might then be somehow transmitted to a receiver and then decoded by
applying $V^*.$ Among all possible left inverses of $V$, we have that
$V^*$ is the unique left inverse that minimizes both the operator norm
and Hilbert-Schmidt norm.

Suppose that in the process of transmission some number, say $m$, of
the components of the vector $Vx$ are lost, garbled or just delayed
for such a long time that one chooses to reconstruct $x$ with what has
been received. In this case we can represent the received vector as
$EVx$, where $E$ is a diagonal matrix of $m$ 0's and $n-m$ 1's
corresponding to the entries of $Vx$ that are, respectively, lost and
received. The 0's in $E$ can be thought of as the coordinates of $Vx$
that have been ``erased" in the language of \cite{GKK}.

There are now two methods by which one could attempt to reconstruct
$x$. Either one is forced to compute a left inverse for $EV$ or one
can continue to use the left inverse $V^*$ for $V$ and accept that $x$
has only been approximately reconstructed.

If $EV$ has a left inverse, then the left inverse of minimum norm is
given by $T^{-1}W^*$ where $EV= WT$ is the polar decomposition and $T
= |EV| = (V^*EV)^{1/2}$.  Thus, the minimum norm of a left inverse is
given by $t_{min}^{-1}$ where $t_{min}$ denotes the least eigenvalue
of $T.$

In the second alternative, the error in reconstructing $x$ is 
given by
$$x - V^*EVx = V^*(I - E)Vx = (I-T^2)x = V^*DVx$$
where $D$ is a diagonal matrix of $m$ 1's and $n-m$ 0's. 
Thus, the norm of the error operator is $1- t_{min}^2$.

Hence we see that, when a left inverse exists, the problems of
minimizing the norm of a left inverse over all frames and of
minimizing the norm of the error operator over all frames are really
equivalent and are both achieved by maximizing the minimal eigenvalue
of $T$.

It is important to note that a left inverse will exist if and only if
the norm of the error operator $V^*DV$ is strictly less than $1$.

In this paper, we study the norms of error operators rather
than those of the left inverses, since this
seems to lead to cleaner formulas and attempt to describe the
frames for which the norms of these error operators are in some sense
minimized, independent of which erasures occur. Of course there are
many ways that one could define ``best'' in this setting and we only
pursue a few reasonable possibilities.

The first quantity that we introduce represents the maximal norm of an
error operator given that some set of $m$ erasures occurs and the
second represents an $\ell^p$-average of the norm of the error operator over the set
of all possible $m$ erasures.

\begin{defn}
  We let $\cl D_m$ denote the set of diagonal matrices that have
  exactly $m$ diagonal entries equal to $1$ and $n-m$ entries equal to
  $0$.
  
  Given an $(n,k)$-frame $F$, we set
  $$e_m^{\infty}(F) = \max\{ \| V^*DV \| : D \in \cl D_m \},$$
  and for
  $1 \le p,$
  $$e_m^p(F)= \{ \begin{pmatrix} n \\ m \end{pmatrix}^{-1} \sum_{D \in
    \cl D_m} \|V^*DV\|^p \}^{1/p}\, ,$$
  where $V$ is the analysis operator
  of $F$ and by the norm of a matrix we always mean its operator norm.
\end{defn}

\begin{remark}
{\rm Note that for a given frame $F$, a left inverse will exist for all
possible sets of $m$ erasures if and only if $e_m^{\infty}(F) < 1.$
Moreover, by the explanation preceding the definition of $e^p_m(F)$, 
whenever at most $m$ erasures occur, then a left inverse of $EV$, 
$L$, will exist satisfying $\|L\| \le
\frac{1}{\sqrt{1- e_m^{\infty}(F)}}.$

In \cite{HP} only the quantity $e_m^{\infty}(F)$ is considered and it is denoted $d_m(F).$

Finally, we remark that the above quantities are invariant under the frame equivalence 
defined in the first section.}
\end{remark}

\begin{defn}
Since $\cl F(n,k)$ is a compact set for any $p$, $1 \le p \le \infty$,
the value 
$$e_1^p(n,k) = \inf \{ e_1^p(F) : F \in \cl F(n,k) \}$$
is attained and we define the {\it (1,p)-erasure frames} to be the
nonempty compact set $\cl E_1^p(n,k)$ of frames where this infimum is
attained, i.e.,
$$\cl E_1^p(n,k) = \{ F \in \cl F(n,k) : e_1^p(F) = e_1^p(n,k) \}.$$
Proceeding inductively, we now set, for $1 \le m \le n$,
$$e_m^p(n,k) = \inf \{ e_m^p(F) : F \in \cl E_{m-1}^p(n,k) \}$$
and define the {\it (m,p)-erasure frames} to be the nonempty compact
subset $\cl E_m^p(n,k)$ of $\cl E_{m-1}^p(n,k)$ where this infimum is
attained.
\end{defn}

In this fashion, we obtain a decreasing family of frames and we 
wish to describe and 
construct the frames in these sets. Since these sets are invariant under frame 
equivalence, we are really only interested in finding representatives for each 
frame equivalence class.

The results of \cite{CK} can be interpreted as characterizing 
$\cl E_1^{\infty}(n,k).$ The following extends their result slightly.

\begin{proposition} \label{prop21}
For $1 < p \le \infty$, the set $\cl E_1^p(n,k)$ coincides with the family of uniform 
$(n,k)$-frames, while $\cl E_1^1(n,k)= \cl F(n,k)$. 
Consequently, for $1 \le p \le \infty, e_1^p(n,k) = k/n.$
\end{proposition}

\begin{proof}

Given an $(n,k)$-frame $F = \{f_1, \ldots ,f_n \},$ if we 
regard the frame vectors as column vectors, then the analysis 
operator $V$ is just the matrix whose $j$-th row is $f^*_j.$ 
Given $D$ 
in $\cl D_1$ which is 1 in the $j$-th entry, we have that
$$\|V^*DV \|^p = \|DVV^*D \|^p = \|f_jf^*_j\|^p = \| 
f_j \|^{2p}.$$

Thus, we see that
$$e_1^{\infty}(F) = \max \{ \|f_j\|^2 : 1 \le j \le n \}$$
and
$$e_1^p(F)= \Bigl(\frac 1 n \sum_j \|f_j\|^{2p}\Bigr)^{1/p}.$$

Since $\sum_j \| f_j \|^2 = tr(VV^*) = k,$ we see that these quantities are minimized when 
$\|f_j\|^2 = k/n$ for all $j$.

Note that when $p=1$, the quantity $e_1^1(F)=k/n$ for every
$(n,k)$-frame and so the result follows.
\end{proof}

We now turn to finding frames that belong to $\cl
E_2^p(n,k).$ By Proposition~\ref{prop21}, for $1 < p$ these are the
uniform $(n,k)$-frames which achieve the infimum of $e_2^p(n,k),$
while for $p=1$ these are just the $(n,k)$-frames that minimize
$e_2^p(n,k)$.

If $D$ is in $\cl D_2$ and has a 1 in the $i$-th and $j$-th diagonal 
entries and $V$ is the analysis operator 
for a uniform
$(n,k)$-frame $F = \{f_1, \ldots ,f_n \}$, then $\|V^*DV\| = 
\|DVV^*D\| = k/n +|\langle f_i,f_j \rangle|
= (1 +\cos(\theta_{i,j}))k/n$ where $0 \le \theta_{i,j} \le \pi/2$ is the 
angle between the lines spanned by the $i$-th and $j$-th frame
vector.

Thus, as observed in \cite{HP} the frames in $\cl E_2^{\infty}(n,k)$
are precisely the uniform $(n,k)$-frames for which the smallest angle
between the lines generated by the frame vectors is maximized. In
\cite{St} these frames were introduced for other reasons and were
called {\it Grassmannian frames.}

We now turn our attention to the frames that are the main topic of this paper.

\begin{defn}
  We call $F$ a {\bf 2-uniform (n,k)-frame} provided that $F$ is a
  uniform $(n,k)$-frame and in addition $\|V^*DV\|$ is a constant for
  all $D$ in $\cl D_2.$
\end{defn}

Unlike uniform frames, 2-uniform frames do not exist for all values of
$k$ and $n$. Later we will give a table that gives a complete list of
all pairs $(n,k)$ for $n \le 50$ for which 2-uniform $(n,k)$-frames
exist over the reals, together with what is known about the numbers of
frame equivalence classes. Each such frame is also a 2-uniform frame
over the complex field, but a complete list of all 2-uniform frames
over the complex field for $n \le 50$ is still not known.

In \cite{HP} it is proved that a uniform $(n,k)$-frame $F$ is
2-uniform if and only if $| \langle f_j,f_i \rangle | = c_{n,k}$ is
constant for all $i \ne j,$ where $$c_{n,k} =
\sqrt{\frac{k(n-k)}{n^2(n-1)}}.$$
The families of frames with this property have been studied 
independently in \cite{St}, where they are called {\it equiangular frames.}

In \cite{HP} it is shown that if there exists a 2-uniform $(n,k)$-frame, 
then every frame in $\cl E_m^{\infty}(n,k)$ is 2-uniform for $2 \le m$ 
and $e_2^{\infty}(n,k) = k/n +c_{n,k}.$ 
If there does not exist a 2-uniform $(n,k)$-frame, then 
necessarily $e_2^{\infty}(n,k) > k/n + c_{n,k}.$

We now prove an analogous result for sufficiently large $p$.

\begin{theorem} \label{thm:largep}
  
  If a 2-uniform $(n,k)$-frame $F$ exists among the uniform frames, 
  then for $p> 2 + \sqrt{\frac{5 k (n-1)}{n-k}}$ and $m\ge 2$, every frame 
  in $\cl E_m^p(n,k)$ is 2-uniform and $e_2^p(n,k) = k/n +c_{n,k}.$ 
  If there does not exist a 2-uniform $(n,k)$-frame, then  
  $e_2^p(n,k) > k/n +c_{n,k}$ for the above values of $p$.

\end{theorem}

\begin{proof} 
  We recall that
  by definition, a frame $F$ in  $\cl E_m^p(n,k)\subset \cl E_2^p(n,k)$ minimizes
  $$2 \sum_{D \in \cl D_2} \|V^*DV\|^p = \sum_{i \ne j}
  (k/n +|p_{i,j}|)^p$$ among all $(n,k)$-frames in $\cl E_1^p(n,k)$,
  where $p_{i,j} = \langle f_i, f_j \rangle$.
  Any such frame $F\in\cl E_1^p(n,k)$ satisfies the 
  constraint 
  $$ 
     \sum_{i \ne j} |p_{i,j}|^2= \frac{nk-k^2}{n}\, ,
  $$   
  because $F$ is uniform and therefore $k = tr(P) = tr(P^2)$.
  To simpify notation, we renumber the $N:=n(n-1)$ 
  quantities $\{|p_{i,j}|,i\neq j\}$ and denote them by $x_l$, $l\in\{1,2, \dots N\}$.
  In addition, we abbreviate $a:=k/n$ and $b:=\sqrt{(nk-k^2)/n}$.
  
  Our proof is a variational argument. It is complete if we show that the minimum of 
  the function $\sigma_p(x):=\sum_{l=1}^N (a+x_l)^p$ subject to $x_l \ge 0$ for all $l$ and
  $\sum_{l=1}^N x_l^2=b^2$ is attained if and only if all $x_l$ are identical. 
  
  As a first step we show that there is $d>0$ such that for any $l$, either $x_l=d$
  or $x_l=0$. 
  For $N=n=2$, this is an exercise in calculus. 
  We want to minimize
  the function $\sigma_p(u,v)=(a+u)^p+(a+v)^p$ subject to the constraints $u^2+v^2=b^2$,
  $u\ge 0$ and $v\ge 0$.
  Our claim is that the minimum occurs either when $u=0$ or $v=0$ or at $u=v$.
  At first we identify critical points of $\sigma_p$ on the arc 
  $A_b=\{u^2+v^2=b^2, u> 0, v> 0\}$. 
  By symmetry, the center $u=v=b/\sqrt 2$ is a critical point and there must be an odd
  number of such points.
  The usual Lagrange multiplier argument shows that at a critical point,
  the equation $(a+u)^{p-1}/u=(a+v)^{p-1}/v$ holds. The set of points
  satisfying this equation in the domain $u,v>0$
  can be split in the three curves $C_1=\{u=v> 0\}$, 
  $C_2=\{v>u> 0, v (a+u)^{p-1} = u(a+v)^{p-1}\}$, and
  $C_3=\{u>v> 0, v (a+u)^{p-1} = u(a+v)^{p-1}\}$. 
  
  We parametrize the curve $C_2$ by $\lambda=(a+u)^{p-1}/u=(a+v)^{p-1}/v$
  and show that this curve intersects only in one point with the arc 
  $A_b$. Once we have proved this, we know this
  critical point cannot be a local minimum, because $\sigma_p(u,\sqrt{b^2-u^2})$ is 
  increasing for sufficiently small values of $u$. By the same argument,  $C_3$ 
  does not contain any local minima, and therefore they must occur either at $u=0$,
  $v=0$, or $u=v$.
   
  To finish the argument for $N=2$, let us assume that $C_2$ and $A_b$ intersect in more
  than one point and derive a contradiction.
  The function $f(t)=\frac{(a+t)^{p-1}}{t}$ has a global minimum at $t=\frac{a}{p-2}$
  so assuming $f(u)=f(v)=\lambda$, $v>u>0$, and $u^2+v^2=b^2$ implies the bounds
  $u<\frac{a}{p-2}$ and $v>\sqrt{b^2-\frac{a^2}{(p-2)^2}}$ for all these intersection points.
  By the condition on $p$, $\frac{a}{p-2}<\frac{b}{\sqrt 2}$, and therefore points in
  $C_2$ sufficiently close to 
  $u=v=\frac{a}{p-2}$ are in the region bounded by the arc $A_b$ and the coordinate axes.
  Moreover, $C_2$ contains points outside of this region, because as the parameter 
  $\lambda$ tends to infinity, so does $v$. 
  All these facts are unchanged if we replace the radius $b$ of the arc by a sufficiently 
  close radius $b'$. Choosing $b'$ appropriately, we can obtain intersection points 
  of $C_2$ and $A_{b'}$ with coordinates $u_2<u_1<v_1<v_2$ such that
  at $(u_1,v_1)$, the radius $\sqrt{u^2+v^2}$ increases as the parameter $\lambda$ increases,
  and at $(u_2,v_2)$ the radius decreases with $\lambda$.
  Consequently,
  $
    \frac{d}{d\lambda} [u^2+v^2]|_{u=u_1,v=v_1} \ge 0
  $
  and
  $
    \frac{d}{d\lambda} [u^2+v^2]|_{u=u_2,v=v_2} \le 0 \, .
  $
  To derive the contradiction, we note that
  \begin{eqnarray*}
     \frac{d}{d\lambda} [u^2+v^2] &=& 2 \left[ \frac{u}{f'(u)}+ \frac{v}{f'(v)} \right]  
%                                   = 2 \left[ \frac{u^3}{(a+u)^{p-2}((p-2)u-a)}+
%                                              \frac{y^3}{(a+y)^{p-2}((p-2)y-a)}\right]
\\
                                  &=& \frac{2}{\lambda^{(p-2)/(p-1)}}
                                      \left[ \frac{u^{2+\frac{2}{p-1} }}{(p-2)u-a}
                                            + \frac{v^{2+\frac{2}{p-1} }}{(p-2)v-a}\right]
  \end{eqnarray*}
  and that
  $g(t):=t^{2+\frac{1}{p-1}}/((p-2)t-a)$
  is negative and strictly decreasing for $t<\frac{a}{p-2}$
  and positive and strictly increasing for $t>(2-\frac 1 p)\frac{a}{p-2}$.
  The condition $p-2>\sqrt{5N}\frac a b$ in our theorem implies 
  that we can choose $b'>\frac{\sqrt{5}a}{p-2}$ and that
  $v_2>v_1>\sqrt{(b')^2-\frac{a^2}{(p-2)^2}}>2\frac{a}{p-2}$. 
  By $g(v_1)+g(u_1)\ge 0$ and the strict monotonicity of $g$  
  we can then infer $g(v_2)+g(u_2)>0$. This contradicts that the radius at $(u_2,v_2)$
  decreases with $\lambda$. Consequently, a second intersection point cannot occur.
  This concludes the first step for the two-variable case.

  The case $N>2$ follows along similar lines. We want to show that all nonzero 
  entries of a minimizing $x$ must be equal. If there are two such non-identical
  entries, say $x_1$ and $x_2$, then one has to be strictly 
  greater than the $\ell^2$-average $b/\sqrt{N}$ of all entries and one has to be 
  strictly less, and thus both 
  add to $x_1^2+x_2^2={\tilde b}^2/N \ge b^2/N$.  Now repeating the argument 
  for $N=2$ with these two entries instead of $u$ and $v$, and $b$ replaced with 
  $\tilde b/\sqrt N$, we see that $\sigma_p(x)$ can be lowered by varying $x_1$ and $x_2$ while
  keeping all other entries fixed, supposing $p-2>\sqrt{5N}\frac{a}{b}$ as stated in the 
  theorem. This contradicts that $x$ is a minimizer.
  
  As the last step we want to ensure that for any $N$, all entries $x_l$ 
  of a minimizer $x$ for $\sigma_p$ are nonzero. Let $M$ denote the number
  of nonzero entries of $x$. By the constraint, we then have 
  $d=b/\sqrt{M}$, and consequently $\sigma_p(x)=M(a+b/\sqrt{M})^p+(N-M)a^p$.
  Assuming $p-2\ge 2 a \sqrt{N}/b=2+2\sqrt{(n-1)k/(n-k)}$ is enough to 
  bound the derivative
  $\partial_M \sigma_p < (a+b/\sqrt{N})^{p-1}(a+b(1-p/2)/\sqrt{N})\leq 0$
  in the interval $1\leq M \leq N$, so
  reducing the number of zeros in $x$ strictly decreases $\sigma_p(x)$. 
  Thus, any minimizer $x$ necessarily satisfies $N=M$ and is therefore unique.
\end{proof}

In an earlier version of this manuscript we erroneously claimed that the above 
result holds for $p> 1$. We are grateful to Srdjan Petrovic for drawing our 
attention to the mistake in our argument. Indeed, for $p=2$ we have 
examples where 2-uniform frames perform the worst among uniform frames. 
We prepare this result with a lemma.

\begin{lemma} \label{prop35}

An $(n,k)$-frame $F = \{f_1, \ldots ,f_n \}$,
is in $\cl E_2^2(n,k)$ if and only if 
it is a uniform frame and minimizes $\sum_{i \ne j} |\langle f_i, f_j \rangle|.$

\end{lemma}

\begin{proof}
  The uniformity of $F$ results from the inclusion $\cl E_2^2(n,k)\subset \cl E_1^2(n,k)$.

  Due to the constraint implicit in the uniformity as stated in the proof of the preceding theorem, 
  $e_2^2(F) = (A +B \sum_{i \ne j} |\langle
  f_i,f_j \rangle|)^{1/2}$ for some positive constants $A$ and $B$
  that depend only on $n$ and $k$. Thus, $e_2^2(F)$ is clearly
  minimized when $\sum_{i \ne j} |\langle f_i,f_j \rangle |$ is
  minimized.
\end{proof} 

\begin{proposition}
If a 2-uniform $(n,k)$-frame $F$ exists among the uniform frames, then it maximizes the 
error $e_2^2(F)$. If a uniform $(n,k)$-frame exists that is not 2-uniform,
then $\cl E_2^2(n,k)$ does not contain any 2-uniform frames.
\end{proposition}
\begin{proof}
This is again a consequence of the fact that 
all frames in $\cl E_2^2(n,k)$ are uniform,
and thus the entries of the associated Grammians observe a constraint 
of the form $\sum_{i\neq j} |p_{i,j}|^2 = b^2$
with a constant $b$ that only depends on $n$ and $k$.

If 2-uniform frames exist, they are maximizers because, subject to the constraint 
$\sum_{l} x_l^2 = b^2$, the function $\sum_{l} |x_l|$  is maximized when all $|x_l|$ are equal.
Thus, given any uniform $(n,k)$-frame that is not 2-uniform, it will perform better and
$\cl E_2^2(n,k)$ cannot contain 2-uniform frames.
\end{proof}

\begin{example}
By Example~\ref{ex:conf} below, we know that 2-uniform $(2k,k)$-frames exist for
infinitely many choices of $k\ge 2$. 
One example of a uniform frame $F'$ that outperforms any such 2-uniform frame is given by basis repetition. 
That is, we simply repeat the vectors of an orthonormal basis twice and rescale them by $1/\sqrt 2$ 
to construct the uniform frame $F'$. For $k\ge 2$, this is not 2-uniform,
because the associated Grammian $P'$ has off-diagonal elements that are zero.
By the preceding proposition then $\cl E_2^2(2k,k)$ does not contain any 2-uniform frames.
\end{example}

\section{Two-uniform frames and graphs} 

In this section we study the existence and construction of 2-uniform
frames.  For many possible values of $(n,k)$, there do not exist any
2-uniform frames. Moreover, when there do exist 2-uniform frames, then
there are at most finitely many such frame equivalence classes and
hence the problem of determining optimal frames in our sense, i.e.,
frames in $\cl E_m^p(n,k)$, is reduced to the problem of finding
representatives for each equivalence class and determining which one
of these finitely many equivalence classes is optimal.

Thanks to the discovery by \cite{St} of the connection between
equiangular frames and the earlier work of Seidel and
his collaborators in graph theory, much of the work of on existence, construction and
determining frame equivalence classes for these frames is already
known and exists in the literature.

We begin this section by summarizing this information. 

\begin{defn}
Given a 2-uniform $(n,k)$-frame $F = \{f_1, \ldots ,f_n \}$ the 
correlation 
matrix is a self-adjoint rank $k$ projection that can be written 
in the form $P = VV^* 
= aI +c_{n,k}Q$ where $a= k/n$, $c_{n,k}$ is given by the formula 
derived in the last 
section, and $Q= (q_{i,j})$ is a self-adjoint matrix  satisfying 
$q_{i,i} = 0$ for all 
$i$ and for $i \ne j, |q_{i,j}| = 1.$ We call 
the $n \times n$ 
self-adjoint matrix $Q$ obtained above the {\bf signature 
matrix} of $F.$

\end{defn}

We shall derive further properties that the signature matrix $Q$ must
satisfy and then use solutions of these equations to generate
2-uniform frames.  The fact that in the real case, $Q$ must be a matrix
of 0's,1's and -1's shows that for fixed $(n,k)$ there are only
finitely many possibilities for the Grammian matrix of a 2-uniform
$(n,k)$-frame. Consequently, up to equivalence, there can be only
finitely many 2-uniform $(n,k)$-frames for each pair $(n,k).$

The key facts about signature matrices are summarized in the following
theorem from \cite{HP}.

\begin{theorem}\cite{HP} \label{thm48}
Let $Q$ be a self-adjoint $n \times n$ matrix $Q$
with $q_{i,i} =0$ for 
all $i$ and $|q_{i,j}| =1$ for all $i \ne j$.  Then the 
following are equivalent:

\begin{itemize}

\item[i)] $Q$ is the signature matrix of a 2-uniform 
$(n,k)$-frame,

\item[ii)] $$Q^2 = 
(n-1)I +\mu Q$$
for some necessarily real number $\mu,$

\item[iii)] $Q$ has exactly two distinct eigenvalues, denoted as $\rho_1 > 
\rho_2.$

\end{itemize}

\end{theorem}

When any of these equivalent conditions hold, then the 
parameters $k$, $\mu$, $\rho_1$, and $\rho_2$ are related by the following equations,
 $$\mu = 
(n-2k)\sqrt{\frac{n-1}{k(n-k)}}
 = \rho_1 +
\rho_2,$$
$$ 1- \rho_1 \rho_2 = n,$$
$$k = n/2 - \frac{\mu n/2}{\sqrt{4(n-1) +
\mu^2}} = \frac{-n\rho_2}{\rho_1 - \rho_2} = mult(\rho_1),$$
where $mult(\rho_1)$ indicates the multiplicity of the eigenvalue $\rho_1$.
In particular, solutions of these equations can only exist for 
real numbers $\mu$ such 
that the formula for $k$ yields an integer.
Moreover, in the case of real 2-uniform frames, the entries of $Q^2$
will all be integers and hence $\mu$ must also be an integer.

The above theorem reduces the construction of 2-uniform 
frames to producing matrices $Q$ satisfying the appropriate equations.

\begin{example} {\it The dimension and codimension 1 case.}  

Let $J_n$ denote the $n \times n$ matrix all of whose entries are 
1. Then $Q = J_n-I_n$ 
satisfies $Q^2 = J_n^2 -2J_n +I_n = (n-2)J_n +I_n = (n-1)I_n +(n-2)Q$ and 
so by our above formulas 
$\mu = n-2$ and $k=1$, yielding the rather uninteresting 
2-uniform frame for $\bb F^1$.

However, $Q = I_n-J_n$ is also a signature matrix with $\mu = 
2-n$, $k=n-1$, which shows that for each $k$ there exists a 
2-uniform $(k+1,k)$-frame.

This frame is described in detail in \cite{CK} and is in fact the only
real uniform $(k+1,k)$-frame, up to some natural equivalence.  We
shall refer to these examples, which exist for every $n$ as the {\it
  trivial 2-uniform frames.}
\end{example}

\begin{example} {\it Conference Matrices. } \label{ex:conf}

The idea of using conference matrices to construct frames of 
this type originates in 
\cite{St}.

A real $n \times n$ matrix $C$ with $c_{i,i} = 0$ and $c_{i,j}= \pm 1$
for $i \ne j$ is called a {\it conference matrix} \cite{CS} provided
$C^2 = (n-1)I.$

Thus, every symmetric conference matrix is a signature matrix 
with $\mu = 0$ and $k = 
n/2$. So, in particular such matrices must be of even size and 
they yield real 
2-uniform $(2k,k)$-frames, for certain values of $k.$

Conference matrices are known to exist for many values of $n$.
Paley\cite{Paley} constructs symmetric
conference matrices, for every $n =
p^j +1 \equiv 2$(mod 4) with $p$ prime. For further examples, see \cite{GS}.

If $C = -C^t$ is a skew-symmetric conference matrix, then 
setting $Q = iC$ yields a 
complex 2-uniform $(2k,k)$-frame.
Similarly, examples of skew-symmetric conference matrices can be found
in many places in the literature. See, for example \cite{Wi} and \cite{GS}.
Note that conference matrices yield 2-uniform frames with redundancy 2.
Conversely, it is not hard to see that the signature matrix of any real 2-uniform 
frame of redundancy 2 is a conference matrix.
\end{example}

\begin{example} 
{\it Hadamard Matrices. }

Using Hadamard matrices to construct 2-uniform frames 
has been discussed in \cite{HP}.
A real $n \times n$ matrix $H$ is called a {\it Hadamard matrix}
\cite{CS} provided that $h_{i,j} = \pm 1$ and $H^*H = nI$. If $H =
H^*$ is a symmetric Hadamard matrix and in addition, $h_{i,i} = 1$ for
all $i$, then $H$ is called a {\it graph Hadamard.}  In this case $Q =
H-I$ is a signature matrix for a real 2-uniform frame with $\mu = -2$
and $k = \frac{n+\sqrt{n}}{2}.$

Similarly, $Q = I-H$ is a signature matrix for a real 2-uniform frame with $\mu = 2$ 
and $k_1 = n-k = \frac{n - \sqrt{n}}{2}$.

Graph Hadamards are known to exist for many values of $n$.
Given two graph Hadamard matrices, their Kronecker tensor product 
gives rise to another graph Hadamard matrix.
Thus, starting with the easily constructed $4 \times 4$ graph
Hadamard, one obtains graph Hadamards of order $4^m$ for every $m$.

A Hadamard matrix $H$ is called a {\it skew Hadamard} if $H +H^* =
2I.$ Note that such a matrix is not actually skew, but is as nearly
skew as a Hadamard matrix can be.

If $H$ is a skew Hadamard matrix, then $Q = \pm i (H-I)$ are signature
matrices for complex 2-uniform frames with $k= n/2.$

Skew Hadamards are known to exist for all $n \equiv 0$(mod 4) and
$n\le 100$ with the exception of $n=72.$

Note that for the 2-uniform frames derived from graph Hadamards, as
$n$ tends to infinity, the redundancy tends to 2, while for all the
skew Hadamards the redundancy is equal to 2.
\end{example}

For further examples of 2-uniform frames, we need to turn to some results in graph theory that were first introduced to the frame theory community by \cite{St}.

\begin{defn}
Given a graph $G$ on $n$ vertices, the {\bf Seidel adjacency 
matrix} of G is defined to be the $n \times n$ matrix $A=  
(a_{i,j})$ where $a_{i,j}$ is defined to be $-1$ when $i$ and $j$ 
are adjacent, $+1$ when $i$ and $j$ are not adjacent, and 0 when 
$i=j$.

Two graphs on $n$ vertices are called {\bf switching 
equivalent} exactly when their Seidel adjacency matrices are 
equivalent via conjugation by an orthogonal matrix that is the product of a 
permutation and a diagonal matrix of $\pm 1$'s.
\end{defn}

The following result, summarized from the results in \cite{HP},
explains the significance of this connection.

\begin{theorem}\cite{HP}
  An $n \times n$ matrix $Q$ is the signature matrix of a real
  2-uniform $(n,k)$-frame if and only if it is the Seidel adjacency
  matrix of a graph with 2 eigenvalues and in this case, $k$ is the
  multiplicity of the largest eigenvalue.  Moreover, if $\{F_i\}$, 
  $i \in\{1,2\}$, is a set of
  real 2-uniform frames, with associated signature matrices $\{Q_i\}$ and
  the corresponding graphs $\{G_i\}$, then $F_1$ and $F_2$ are frame
  equivalent if and only if $G_1$ and $G_2$ are switching equivalent
  graphs.
\end{theorem}

There is a considerable literature in graph theory dedicated to
finding graphs with two eigenvalues and classifying these graphs up to
switching equivalence. By referring to this literature, we can give a
complete list of all integers $n \le 50$ for which such graphs (and
hence 2-uniform frames) are known to exist, together with information
about how many frame equivalence classes there are in each case.

This information is gathered together in {\bf Table 1.} When an
integer $j$ appears in the column labeled, ``frame equivalence
classes'', it indicates that exactly $j$ inequivalent real 2-uniform
$(n,k)$-frames exist. When the symbol $j+$ appears, it indicates
that at least $j$ inequivalent real 2-uniform $(n,k)$-frames are known
to exist, but it is not known yet if this is exhausts all equivalence
classes. The letters $C,H$ and $G$ in the column labeled ``type''
indicate that the corresponding frames are all constructed using
conference matrices, graph Hadamards, or only arise from certain
graphs, respectively.

So for example, using Table 1, and looking at $n=36,$ we see that
there exist at least 227 switching inequivalent graph Hadamard
matrices and these can be used to construct at least 227 frame
inequivalent 2-uniform $(36,15)$-frames and at least 227 frame
inequivalent 2-uniform $(36,21)$-frames.  For $n=276$, there
exists a graph whose Seidel adjacency matrix has exactly 2
eigenvalues, but it is neither a conference matrix nor graph Hadamard
matrix, and this matrix can be used to construct a 2-uniform
$(276,23)$-frame that up to frame equivalence is the unique such
frame.

\begin{table}[bth]\label{eterms}

\caption{Real 2-Uniform Frames}

  \begin{center}

\begin{tabular}
 {|c|c|c|c|}\hline
No. of vectors $n$ &Dimension $k$  & No. of Equiv. Classes & Type \\
\hline

6 & 3 & 1 & C \\

10 & 5 & 1 & C \\

14 & 7 & 1 & C \\

16 & 6;10 & 1 & H\\

18 & 9 & 1 & C \\

26 & 13 & 4 & C \\

28 & 7;21 & 1 & G\\

30 & 15 & 6+& C\\

36 & 15;21 & 227+& H\\

38 & 19 & 11+& C\\

42 & 21 & 18+& C\\

46 & 23 & 80+& C\\

50 & 25 & 18+& C\\

176 & 22;154 & 1 & G\\

276 & 23;253 & 1 & G\\ \hline

\end{tabular}

  \end{center}

\end{table}

The number of equivalence classes are often computed by using the
theory and enumeration of two-graphs. A {\it two-graph} 
$(\Omega, \Delta)$ is a pair consisting of a 
vertex set $\Omega$ and a collection $\Delta$ of three element 
subsets of $\Omega$ such that every four element subset of 
$\Omega$ contains an even number of the sets from $\Delta.$
A two-graph is regular, provided that every two element subset 
of $\Omega$ is contained in the same number, $\alpha$, of sets 
in $\Delta.$

Given $n$, Seidel~\cite{S} exhibits a one-to-one correspondence 
between the two-graphs on the set of $n$ elements and the 
switching equivalence classes of graphs on $n$ elements and 
gives a concrete means, given the two-graph, to construct a 
graph from the corresponding switching class. 
Thus, a two-graph can be regarded as a switching equivalence 
class of ordinary graphs.

In \cite{St}, it was noted that signature matrices of real 
2-uniform frames are always Seidel adjacency matrices of 
regular two-graphs. The following result from \cite{HP} more fully 
summarizes this connection.

\begin{theorem}\cite{HP} An $n \times n$ matrix $Q$ is the signature matrix 
of a real 2-uniform $(n,k)$-frame if and only if it is the 
Seidel adjacency matrix of a graph on $n$ vertices whose 
switching equivalence class is a regular two-graph on $n$ 
vertices
with parameter $\alpha$. This relationship defines a one-to-one 
correspondence between frame equivalence classes of real 
2-uniform frames and regular two-graphs.

\end{theorem}

The relationship between the parameter $\alpha$ and earlier 
parameters is given by the equations, $$-2\alpha = (1+
\rho_1)(1+\rho_2) = 2+\mu +n.$$

Thus, by the above theorem every regular two-graph produces a 
real 2-uniform frame. For a given $n$ these could just be the 
trivial, known examples corresponding to $k=n-1,1.$
In \cite{S} many of the known regular two-graphs are listed and 
it is elementary to use the formulas given above to determine 
the pairs $(n,k)$ for which they yield a real 
2-uniform frame.

In particular, the two-graph $\Omega^{-}(6,2)$ yields a 
2-uniform $(28,7)$-frame and a 2-uniform $(28,21)$-frame. Since by Table 1 there is only one such
equivalence class, the frames derived from this two-graph will be
frame equivalent to the frames generated by the $28 \times 28$ signature 
matrix (and
its negative)
found by Holmes~\cite{Ho}.

\begin{prob} One fundamental question that we have not been able to answer 
  is whether or not regular two-graphs exist that give rise to 2-uniform 
  frames of arbitrarily large redundancy. The examples that come from 
  conference matrices all have redundancy 2 and those coming from Hadamard 
  matrices have redundancy approaching 2. The existence of two-graphs
  with arbitrarily large redundancy may possibly be a consequence of Ramsey 
  theory.
\end{prob}

\section{Graphs and Error Bounds}

In this section we derive estimates and formulas for $e_m^p(F)$
when $F$ is a real 2-uniform frame, using connectivity properties of the
graph associated to the signature matrix $Q$ of $F.$

Recall that if $F$ is a real 2-uniform $(n,k)$-frame, then the Grammian $P$ of
$F$ is an $n \times n$ matrix that is a projection of rank $k$ and has
the form
$P= k/nI +c_{n,k}Q$, where $c_{n,k}= \sqrt{ \frac{k(n-k)}{n^2(n-1)}}$ and
$Q$ is the Seidel adjacency matrix of a certain graph that we shall
denote $G_F.$

We also need to recall a few concepts from graph theory. 

\begin{defn}
A graph on
$m$ vertices is called  {\bf complete bipartite} provided that the
vertex set can be divided into two disjoint subsets, $V_1$ and $V_2$,
of sizes, say $m_1$ and $m_2$ with $m_1 +m_2 = m$, such that no pair
of vertices in $V_1$ or in $V_2$ are neighbors, but every vertex in
$V_1$ is adjacent to every vertex in $V_2$. We shall denote this
graph by $B(m_1,m_2)$.  {\bf In our definition of complete bipartite graph,
we allow the possibility that one of the sets is empty}, so that the
complete bipartite graph, $B(0,m)$, is the graph on $m$ vertices with
no edges.  If $G$ is a graph with vertex set $V$ and $W \subset V$
then by the {\it induced subgraph on W}, we mean the graph with vertex
set $W$ and two vertices in $W$ are adjacent if and only if they were
adjacent in $G$.
\end{defn}

Note that if $F$ is a real 2-uniform $(n,k)$-frame with signature matrix
$Q$ and graph $G_F,$ then the Seidel adjacency matrix of an induced
subgraph on $m$ vertices is just the $m \times m$ matrix obtained by
compressing $Q$ to the corresponding entries.

We are grateful to Ryan Pepper for the following observation, which
can also be found in the work of Seidel.

\begin{lemma}
A graph on $m$ vertices is switching equivalent to the graph with no
edges if and only if it is complete bipartite.
\end{lemma}
\begin{proof}
Given any complete bipartite graph corresponding to a preselected partition
of vertices into two sets, we show that it can be obtained by switching 
the empty graph on $m$ vertices. Without loss  of generality we may order 
the vertices $\{v_1,v_2, \dots v_m\}$
such that the partition is given by $\{v_j\}_{j\leq r}$ and $\{v_j\}_{j>r}$,
with $0\leq r \leq m$. Let us choose the switching matrix
$$
   S=\left(\begin{array}{cc} +I_r & 0 \\ 0 & -I_{m-r} \end{array}\right) \, .
$$
   The empty graph on $m$ vertices has the Seidel adjacency matrix 
$J_m-I_m$, so that of the switched graph is then
$$
  S(J-I)S = \left(\begin{array}{cc} J_r-I_r  &-J_{r,m-r} \\ -J_{m-r,r} & J_{m-r}-I_{m-r} 
\end{array}\right)
$$
which by inspection belongs to the preselected complete bipartite graph.
Moreover, switching the empty graph always leads to complete bipartite graphs.  
Again, the empty (edgless) graph is 
in our sense bipartite, corresponding to a partition $\emptyset$ and 
$\{v_1,v_2,\dots v_m\}$.
\end{proof}

\subsection{Error estimates for $e^\infty_m(F)$.}

\begin{theorem} \label{thm:cb}
  Let $F$ be a real 2-uniform $(n,k)$-frame. Then
  $e_m^{\infty}(F) \le k/n +(m-1)c_{n,k}$ with equality if and only
  if $G_F$ contains an induced subgraph on $m$ vertices that is
  complete bipartite.
\end{theorem}

\begin{proof}
The proof has three parts. First we show that the operator norm $||DVV^*D||$
is equal to the largest eigenvalue of the compression of $VV^*$ to the
rows and columns where $D$ has $1$'s. Then we bound the largest possible eigenvalue.
Finally we characterize the case when this bound is saturated.

To begin with, we note that $DVV^*D$ is a positive operator,
and so is its compression $(VV^*)_m$ to the rows and columns where $D$ has
$1$'s. Consequently, the operator norm of $DVV^*D$ is equal to the largest 
eigenvalue of $(VV^*)_m$.
This eigenvalue, in turn, follows from the largest eigenvalue
of the compression of $Q$, because $(VV^*)_m=\frac{k}{n} I_m + c_{n,k}Q_m$.
So we can reduce the calculation to that of the largest eigenvalue of
$Q_m$. In fact, to simplify the argument, we will look for the largest
possible eigenvalue of $Q_m+I_m$ and later adjust for the added constant.

We now claim that the largest eigenvalue of $Q_m+I_m$ occurs when
$Q_m+I_m=J_m$, that is, the matrix of all $1$'s. This follows from 
considering that for any given $x\in \R^m$ and $Q_m$,
   changing signs to make all their entries positive increases
   $\langle x,(Q_m+I_m) x\rangle/||x||^2$.
   
By inspection, the largest eigenvalue of $J_m$ is $m$, so that 
of $J_m-I_m$ is $m-1$, and the claimed error bound follows:
$$
  ||DVV^*D|| \leq \frac k n  +(m-1)c_{n,k} \, .
$$

To characterize cases of equality, suppose $G$ contains an induced
subgraph that is switching equivalent to the graph with no edges. If
we choose $D$ to have $1$'s in the places on the diagonal
corresponding to the vertices of this subgraph and $0$'s everywhere
else, then $D(I+Q)D$ is switching equivalent to $DJ_nD$ and so the error is
$e_m^{\infty}(F) = k/n +(m-1)c_{n,k}$. Conversely, assume that
equality holds in the error estimate. Then, $||D(I+Q)D||=m$. Given an
eigenvector $x$ corresponding to eigenvalue $\pm m$, we may choose a
switching matrix $S$ such that all of the entries of $Sx$ are
positive. Similarly to the above reasoning, all the entries in
$S(I+Q)S$ must be $1$'s in the rows and columns where $D$ has $1$'s on
the diagonal, otherwise it would be possible to increase the largest
eigenvalue of $SD(I+Q)DS$ by flipping signs in $Q$, contradicting that
the inequality is saturated. Hence, the induced subgraph on these
vertices is switching equivalent to the edgeless graph.
\end{proof}

\begin{corollary}
For a given $m$, a real 2-uniform frame $F$ maximizes the error 
$e_m^\infty(F)$ iff $G_F$ has an induced subgraph 
on $m$ vertices that is compelete bipartite.
\end{corollary}

\begin{corollary} Let $F$ be a real 2-uniform $(n,k)$-frame. 
Then either $G_F$ contains an induced complete bipartite graph on 3 vertices
or it is switching equivalent to the complete graph on $n$ vertices. 
Consequently, if $k< n-1$ we have
  $e_3^{\infty}(F) = k/n +2c_{n,k}$.
\end{corollary}
\begin{proof}
  Let us assume that $G_F$ has no induced complete bipartite subgraph on 3 vertices.
  We may choose one vertex and switch the others if necessary
  in order to have edges between this one and all others.  
  Then any two vertices must be adjacent, otherwise there would be
  an induced complete bipartite subgraph on 3 vertices. Thus, the resulting graph is 
  a switched version of $G_F$ that is complete. This corresponds to $F$ being equivalent 
  to the uniform $(n,n-1)$-frame.
\end{proof}

\begin{corollary}
If $F_1$ and $F_2$ are real 2-uniform $(n,k)$-graphs, then
$e_3^\infty(F_1)$ $=e_3^\infty(F_2)$.
\end{corollary}

By analogy with our earlier definitions, we call a Parseval frame $F$ {\it
  3-uniform} if it is 2-uniform and if the error $\|DVV^*D\|$ associated
with a deletion is constant for every $D \in \cl D_3$. 

\begin{corollary}
The trivial 2-uniform $(n,k)$- frames, corresponding to $k=1$ 
or $k=n-1$, are 3-uniform. Conversely,
if $F$ is a real $3$-uniform $(n,k)$-frame, then either $k=1$ or $k=n-1$
and it is equivalent to the corresponding trivial frame. 
\end{corollary}
\begin{proof}
It is clear from their definition that the trivial frames are 3-uniform.
What we need to show is that if $F$ is 3-uniform then $G_F$ is either
switching equivalent to the complete or to the edgeless graph.

To begin with, we pick a vertex and switch the others if necessary in order
to isolate it. Any two additional vertices are either adjacent or not,
and computing the norm of $DVV^*D$, where $D$ is associated with these 3
vertices, distinguishes these
two cases. However, 3-uniformity then implies that {\it every} additional pair
of vertices must behave the same way. Thus, if there is one edgeless induced subgraph
on 3 vertices, then 
all of $G_F$ is edgeless. On the other hand, if there is one neighboring pair, 
then all pairs of
vertices except those including the isolated one are neighbors. Switching this
one vertex then yields the complete graph. 
\end{proof}

We now discuss how non-existence of complete bipartites gives rise to
refined error bounds.

\begin{defn} 
Let ${\cl G}_m^{(s)}$ denote the set of graphs on $m$
vertices such that $s$ is the minimum number of edges occuring among graphs 
in the equivalence class associated with each $G \in 
{\cl G}_m^{(s)}$.
\end{defn}

We note that ${\cl G}_m^{(0)}$ are the complete bipartites on $m$ vertices,
and ${\cl G}_m^{(1)}$ is the equivalence class of graphs that may be reduced
to one edge by switching. However, for larger values of $s$, 
${\cl G}_m^{(s)}$ may contain more than one equivalence class.

\begin{lemma}
Let, $Q^{(0)}$, $Q^{(1)}$ and $Q^{(s)}$, $s\ge 1$, be Seidel adjacency 
matrices belonging to graphs $G^{(0)} \in {\cl G}_m^{(0)}$, 
$G^{(1)}\in {\cl G}_m^{(1)}$, and $G^{(s)}\in {\cl G}_m^{(s)}$,
respectively, for some common number of vertices $m\ge 3$. 
Denote by $\lambda^{(0)}$, $\lambda^{(1)}$ and 
$\lambda^{(s)}$ the largest eigenvalue of $Q^{(0)}$, $Q^{(1)}$ and 
$Q^{(s)}$. Then $\lambda^{(0)}\ge\lambda^{(1)}\ge\lambda^{(s)}$.
\end{lemma}
\begin{proof}
By appropriate switching, we can make $G^{(0)}$, $G^{(1)}$ and $G^{(s)}$ have a
minimal number of edges in their respective equivalence class. In particular,
then $G^{(0)}$ is the edgeless graph. Permuting the vertices if necessary, we 
have $G^{(0)}\subset G^{(1)} \subset G^{(s)}$. 
To simplify notation, we can choose this permutation in such a way that there is
an edge between the $m$-th and $m-1$-th vertex in $G^{(1)}$ and in $G^{(s)}$.
Since switching
corresponds to a change of basis in ${\mathbb R}^n$, the maximal eigenvalues of 
the Seidel adjacency matrices $Q^{(0)}$, $Q^{(1)}$ and $Q^{(s)}$ are unchanged. 
After switching, the components of $Q^{(1)}$ and $Q^{(s)}$ observe 
$q^{(1)}_{m,m-1}=$ $q^{(1)}_{m-1,m}=$ $q^{(s)}_{m,m-1}=$ $q^{(s)}_{m-1,m}=-1$.
The inequality between the largest eigenvalues of $Q^{(0)}$ and $Q^{(1)}$
follows by explicit computation, $\lambda^{(0)} = m - 1 \ge \lambda^{(1)}=\frac{m}{2}-2+
\sqrt{\frac{m^2}{4}+m-3}$ for $m\ge 3$.

To establish the inequality $\lambda^{(1)}\ge \lambda^{(s)}$, we use a 
variational argument similar to that in the proof of Theorem~\ref{thm:cb}.
We consider a normalized eigenvector $x$ belonging to the largest eigenvalue of 
$Q^{(s)}$. We show there is a normalized vector $p$ such that
$\lambda^{(1)}\ge\langle p, Q^{(1)} p\rangle \ge$
$\langle x, Q^{(s)} x\rangle=\lambda^{(s)}$.
The construction of $p$ proceeds in several parts: \\
{\it Part 1.} First let us assume 
that not all of $\{x_i\}_{i=1}^{m-2}$
are negative, otherwise we multiply $x$ by $-1$. Now we set $x'_i=|x_i|$ for $i\leq m-2$
and replace $q_{i,j}$ by $q'_{i,j}:=|q_{i,j}|$ for the block given by 
$i,j \leq m-2$. 
At the same time, we modify the last two rows and columns of $Q^{(s)}$ as follows.
If $x_{m-1}$ and $x_m$ are both positive or both negative, we set $p_i=|x_i|$ for all $i$,
let $q'_{m-1,j}=q'_{j,m-1}=1$ and
$q'_{m,j}=q'_{j,m}=1$ for $j\leq m-2$, and thus obtain $Q^{(1)}=Q'$ from 
$Q^{(s)}$ while only increasing $\langle x, Q^{(s)} x \rangle$ to $\langle p, 
Q' p\rangle$, which finishes the proof for this case.\\
{\it Part 2.}
If one component of $x$, say $x_m$, is negative, we 
set $x'_i= |x_i|$ only for $i \leq m-1$ and retain $x'_m=x_m$, while letting 
$q'_{i,m}=q'_{m,i}=-1$ if $x_i<0$ and $q'_{i,m}=q'_{m,i}=1$ if $x_i>0$, for $i \leq m-2$. 
This ensures that $\langle x', Q' x'\rangle \ge \langle x, Q^{(s)} x \rangle$,
while we have reduced the occurence of $-1$'s in $Q'$ to the last row and column.
Note that by 
our assumption that 
not all $\{x_i\}_{i=1}^{m-2}$ are negative, according to this procedure, 
there remains at least one entry $q'_{i,m}=q'_{m,i}=1$ in the last row and column of $Q'$.\\ 
{\it Part 3.}
Now define $x''$ by $x''_i=|x_i|$ and let $Q''$ be obtained from switching signs in the last row
and column of $Q'$. Then $\langle x', Q' x'\rangle=\langle x'', Q'' x''\rangle$. 
Since $x''$ has all positive entries and $Q''$ contains at least one pair of $-1$'s in the last 
row and column, setting all entries to $1$ but this one pair $q''_{i,m}=q''_{m,i}=-1$ only increases 
$\langle x'', Q'' x''\rangle$ and transforms $Q''$, together with a permutation of indices,
to $Q^{(1)}$. Applying the same permutation to the entries of $x''$ yields $p$ satisfying 
$\langle p, Q^{(1)} p\rangle \ge \langle x, Q^{(s)} x\rangle$.
\end{proof}

\begin{theorem}\label{cor555}
Given a real 2-uniform  $(n,k)$-frame $F$ such that for some $m\ge 3$, the 
associated graph $G_F$ does not have an induced complete bipartite subgraph
on $m$ vertices. Then we have the improved error bound
$$
  e_m^\infty(F) \leq \frac{k}{n}+c_{n,k}\Bigl(\frac{m}{2}-2+
    \sqrt{\frac{m^2}{4}+m-3}\Bigr)\, .
$$
If $G_F$ contains an induced subgraph on $m$ vertices that differs
from a complete bipartite by one edge, then equality holds. 
\end{theorem}
\begin{proof}
The improved error bound results from the fact that by the preceding lemma,
in the absence of complete bipartites on $m$ vertices, the graphs 
in ${\cl G}_m^{(1)}$ maximize the error.

To prepare the argument, we recall that $VV^*=\frac k n I_n +
  c_{n,k} Q$ is a projection, so the
compression of $DVV^*D$ to the rows and columns where
$D$ has $1$'s, henceforth denoted as $\frac k n I_m+
  c_{n,k} Q_m$,  is a non-negative operator.
Consequently, the norm of $DVV^*D$ equals that of $\frac k n I_m +
  c_{n,k} Q_m $ and is given by its largest eigenvalue.
To obtain this eigenvalue, it is enough to consider $Q_m$.

In the absence of complete bipartites, we know that 
any matrix $Q_m$, obtained from the compression of $Q$ to $m$ rows
and columns, corresponds to a graph $G\in {\cl G}_m^{(s)}$ with $s\ge 1$. 
By the inequality in the preceding lemma, to bound the 
largest possible eigenvalue we only need to consider $s=1$.
Since we may switch and permute $G$ without loss of generality, we again 
choose $Q_m$ to be the $m \times m$ matrix of all $1$'s, except the $0$'s 
on the diagonal and the two elements $q_{m-1,m}=q_{m,m-1}=-1$. By inspection, 
$Q_m$ has an eigenvector $(0,0,\dots,1,-1)$ with eigenvalue $1$, a set 
of $m-3$ linearly independent eigenvectors of eigenvalue $-1$ given 
by $(1,-1,0,\dots,0)$ and those obtained when exchanging its second entry 
with all others except
the first one and the last two entries.
The larger of the remaining two eigenvalues is $\lambda=\frac{m}{2}-2+ 
\sqrt{(\frac{m}{2}-1)^2+2m-4}$ which is seen to be greater than or equal 
to $1$ for
$m\ge 3$. The claimed error bound follows.
\end{proof}

The next five results can be deduced by converting results in
\cite{GS2}, especially Remark~2.8 and Theorem~2.7, into statements
about frames and re-deriving their formulas in terms of the parameter $k$, 
but it is perhaps clearer to deduce them directly. A main part of the 
results consists of sufficient conditions that rule out the existence 
of induced complete bipartite subgraphs on certain numbers of vertices.
Theorem~\ref{thm:557} summarizes these conditions.

\begin{proposition} \label{prop57}
Let $F$ be a real 2-uniform $(n,k)$-frame, and suppose 
$1 +\sqrt{\frac{(n-k)(n-1)}{k}} < m$.
Then the associated graph 
$G_F$ does not contain an induced subgraph on $m$ vertices that is complete bipartite.
\end{proposition}
\begin{proof}
  Since $e_m^{\infty}(F)\leq 1$ and $\frac k n +(m-1)c_{n,k}>1$ for $1
  +\sqrt{\frac{(n-k)(n-1)}{k}} < m$, we have $e_m^{\infty}(F)< \frac
  k n +(m-1)c_{n,k}$. This excludes an induced complete bipartite 
  subgraph, because otherwise equality would hold.
\end{proof}
\begin{corollary} \label{cor:556}
Let $F$ be a real 2-uniform $(n,k)$-frame. If $n-k+1 \leq m \leq n,$ then
no induced subgraph of $G_F$ on $m$ vertices can be complete
bipartite.
\end{corollary}
\begin{proof}
Since $1+\sqrt{\frac{(n-k)(n-1)}{k}} < n-k+1$, there cannot be any 
induced complete bipartite subgraphs of $m$ vertices 
when $n-k+1\leq m \leq n$.
\end{proof}

\begin{proposition} Let $F$ be a real 2-uniform $(n,k)$-frame.
  If $m>k$, then the Seidel adjacency matrix of any 
  induced subgraph of $G_F$ on $m$ vertices has an eigenvalue $-\frac{k}{n c_{n,k}}$.
\end{proposition}
\begin{proof}
  Take $m>k$ and consider
$P=\left( \begin{array}{cc}\frac k n I_m +c_{n,k} Q_m & *\\ * & * \end{array}\right)$,
where $*$ denotes the remaining entries of $P$ outside of the first $m$ rows and columns. 
Since $\mathrm{rk}(P)=k$, these $m$ columns are linearly dependent and
$$
  0 \in \sigma(\frac k n I_m +c_{n,k} Q_m )\, .
$$
Thus
$$
  0=\frac k  n +\lambda c_{n,k}
$$
for some eigenvalue $\lambda$ of $Q_m$.
\end{proof}
As a precursor to the next consequence, we 
recall that since $VV^*$ is an $n\times n$ matrix and a projection of
rank $k$, it has eigenvalues $1$ with multiplicity $k$ and $0$ with
multiplicity $n-k$. 
Hence $Q$ has eigenvalues $\rho_2<\rho_1$ with
$\rho_1=\frac{n-k}{nc_{n,k}}=\sqrt{\frac{(n-k)(n-1)}{k}}$ of multiplicity
$k$ and $\rho_2=-\frac{k}{n c_{n,k}}=-\sqrt{\frac{k(n-1)}{n-k}}$ of
multiplicity $n-k$.
\begin{corollary} \label{cor511}
Given a real 2-uniform $(n,k)$-frame $F$, then $G_F$ has
no induced subgraph on $m>k$ vertices that is complete bipartite.
\end{corollary}
\begin{proof}
If it had, then the signature matrix $Q_m$ associated with the subset
of vertices that form the induced complete bipartite subgraph 
would have eigenvalues $\sigma(Q_m)=\sigma(J_m-I_m)=\{-1,m-1\}$,
so $\frac{k}{n c_{n,k}}=1$, contradicting $k/n>c_{n,k}$.
\end{proof}

\begin{proposition} Let $F$ be a real 2-uniform $(n,k)$-frame.
If $m>n-k$, then the Seidel adjacency matrix of every induced subgraph on $m$ 
vertices has an eigenvalue $\frac{n-k}{nc_{n,k}}$. 
\end{proposition}
\begin{proof}
 The projection onto the complement of the range of $P$, $I-P$, has rank
$n-k$. So $0\in \sigma(I_m-(\frac k n I_m +c_{n,k} Q_m))$
and $0=1 - \frac k n - \lambda c_{n,k}$ for some eigenvalue $\lambda$ of $Q_m$.
\end{proof}
\begin{corollary} \label{cor513} Let $F$ be a real 2-uniform $(n,k)$-frame.
If $m>n-k,$ then no induced subgraph of $G_F$ on $m$ vertices is complete bipartite.
\end{corollary}
\begin{proof}
$m>n-k$ gives $m-1\ge n-k$. But $nc_{n,k}=\sqrt{\frac{k(n-k)}{n-1}}>1$
so $\lambda = \frac{n-k}{nc_{n,k}} < n-k \leq m-1$. Thus $\lambda \not\in \{m,m-1\}$.
\end{proof}

\begin{theorem} \label{thm:557}
Let $F$ be a real 2-uniform $(n,k)$-frame. If $G_F$ contains an
  induced subgraph on $m$ vertices that is complete bipartite then
$m \le \min \{ k,n-k, 1+\sqrt{\frac{(n-k)(n-1)}{k}} \}$
\end{theorem}
\begin{proof}
Follows from preceding Proposition~\ref{prop57}, Corollaries~\ref{cor511}
and \ref{cor513}. 
\end{proof}

The next result shows which number of erasures may cause a 2-uniform frame
to lose all information contained in some encoded vector. 

\begin{proposition}
For $n-k+1 \leq m \leq n$ and any (real or complex) 2-uniform $(n,k)$-frame $F$, 
$e_m^{\infty}(F)=1$.
\end{proposition}
\begin{proof}
This follows from an eigenvalue interlacing theorem and the multiplicity $k$ 
of the eigenvalue one of $P=VV^*$. If $k\ge 2$ and
$m=n-1$, then the $k-1$ largest eigenvalues of $DVV^*D$
must lie between the $k$ largest eigenvalues of $VV^*$, which are all one.
By iteration, the eigenvalue one will remain up to $m=n-k+1$.
\end{proof}
\medskip
\subsection{Computation of the error $e^p_3(F)$.}~

\noindent We now turn our attention to $e_3^p(F)$. Recall
that switching equivalence leads to two different equivalence classes 
$\Gamma_e$ 
and $\Gamma_o$ for 3-vertex graphs, those with an even number and those with an odd 
number of edges, repectively. We observe that $\Gamma_e$ contains exactly the 
complete bipartite graphs with 3 vertices.

\begin{lemma} \label{lem521}
The number of complete bipartite 3-vertex subgraphs $E_3(G)$ in a graph $G$ that 
corresponds to a real 2-uniform $(n,k)$-frame depends only on $n$ and $k$. It is given
by
$$
  E_3(G)=\Bigl({n \atop 3}\Bigr)-\frac{v(n-1)c}{6}- (n-2v+c)\frac{(n-1) v}{2}
$$ 
with 
$$ v = \frac 1 2 (n-2-\sqrt{\frac{(n-k)(n-1)}{k}}+\sqrt{\frac{k(n-1)}{n-k}})$$
and $$c=v-1-\frac 1 4 (\sqrt{\frac{(n-k)(n-1)}{k}}-1)(\sqrt{\frac{k(n-1)}{n-k}}+1). $$ 
\end{lemma}
\begin{proof}
By Seidel, if $G$ is a graph in the switching class of a regular two-graph 
and if $G$ has an isolated vertex, then the induced graph $\tilde G:=G\setminus \{\omega\}$
is strongly regular 
\cite[Theorems 6.11 and 7.2]{S}. Thus, $\tilde G$ is characterized by
the tuple $(n-1,v,p,q)$ which represent, respectively,
the total number of vertices $n-1$, the common
valency $v$ of each vertex, whenever two vertices are adjacent there
are $p$ vertices adjacent to one vertex and not the other, and
whenever two vertices are non-adjacent there are $q$ vertices adjacent
to one vertex and not the other. If we let $c$ denote the number of
common neighboring vertices of two adjacent vertices, then $c+p+1=v.$

As a first step, we count the number of odd-edged induced subgraphs in $\tilde G$,
denoted as $O_3(\tilde G)$. The total number of edges in $\tilde G$ is $\frac{(n-1) v}{2}$. 
Each edge belongs to $c$ 3-edged graphs.
Therefore, $\frac{(n-1)v}{2}c$ counts each 3-edged graph three times
and so there are $\frac{v(n-1)c}{6}$ 3-edged subgraphs.
To arrive at the number of 1-edged subgraphs, we recall that
two connected vertices are each connected to $v-1$ 
other vertices, and have $c$ of these as common neighbors. Thus, these vertices
are connected to $2(v-1)-c$ other vertices, and not connected to $(n-1)-2v+c$. Hence, there
are $(n-1-2v+c)\frac{(n-1)v}{2}$ 1-edged subgraphs. The number of 3-vertex subgraphs in
$\tilde G$ with
odd edges is consequently 
$$
   O_3(\tilde G)=\frac{v (n-1) c}{6} +(n-1-2v+c)\frac{(n-1)v}{2} \, .
$$
The number of odd-edged 3-vertex subgraphs in $G$ is then 
$$
  O_3(G) = O_3(\tilde G) +\frac{(n-1)v}{2}
$$
due to adding an 1-edged subgraph for every edge in $\tilde G$ when including
the isolated vertex.
Thus $O_3(G)$ and $E_3(G)=({n \atop 3})-O_3(G)$ follow once $v$ and $c$ are given.
What remains is to deduce their values from the 2-uniformity of the 
frame belonging to $G$.

By \cite[Theorem 7.5]{S} if we pick $G$ from the switching class of a regular two-graph
such that $G$ has an isolated vertex $\omega$, then the Seidel adjacency matrix of 
$\tilde G=G\setminus\{\omega\}$ 
has eigenvalues $\rho_0=\rho_1+\rho_2$, $\rho_1$, and $\rho_2$, where 
$\rho_1=\sqrt{\frac{(n-k)(n-1)}{k}}$
and $\rho_2=-\sqrt{\frac{k(n-1)}{n-k}}$ are the eigenvalues associated with $G$.
By Seidel (pp.\ 155--156), the valency is $v=\frac 1 2(n-2-\rho_0)$ and
$4p=-(\rho_1-1)(\rho_2-1)$, which concludes the proof.
\end{proof}

\begin{proposition}\label{lem517}
For a real 2-uniform $(n,k)$-frame $F$, the error $e_3^p(F)$ for $2
\le p < \infty$ is given by
$$
 e_3^p(F) = \left( \Bigl({n \atop 3}\Bigr)^{-1} 
        \left[ (\frac k n +2 c_{n,k})^p E_3(G)
        +(\frac k n +c_{n,k})^p O_3(G) \right] \right)^{1/p}\, , 
$$ 
where $E_3(G)$ and $O_3(G)$ are the constants that were calculated in the preceding lemma. 
\end{proposition}
\begin{proof}
For an induced 3-vertex subgraph in $\Gamma_e$ associated 
with a 2-uniform $(n,k)$-frame $F$ and $D\in{\cl D}_3$, an explicit computation gives
$||DVV^*D|| = \frac k n +2 c_{n,k}$, whereas if the subgraph is in $\Gamma_o$ then
$||DVV^*D||=\frac k n +c_{n,k}$. 
Consequently,
%$$
%  e_3^p(F) = \left( \Bigl({n \atop 3}\Bigr)^{-1} 
%        \left[ (\frac k n +2 c_{n,k})^p E_3(G)
%        +(\frac k n +c_{n,k})^p O_3(G) \right] \right)^{1/p}\, , 
%$$
the definition of $e_3^p(F)$ simplifies to the claimed expression.
\end{proof}

\begin{corollary}
For any two real 2-uniform $(n,k)$-frames $F_1$ and $F_2$,
$e_3^p(F_1)$ $=e_3^p(F_2)$ for $2 \le p < \infty$.
\end{corollary}
\begin{proof}
Follows from Lemma~\ref{lem521} and Proposition~\ref{lem517}. 
\end{proof}

\section{Error Estimates for Concrete Frames}

In this section, we use the inequalities and methods of the previous
section to explicitly compute the error estimates, $e_m^p(F)$ for
various 2-uniform $(n,k)$-frames. In addition, we investigate how
the error estimates compare to an explicit, computer-aided calculation 
of the error.

%\subsection{Conference matrices}
We begin with an example of frames constructed with the help
of conference matrices.

\begin{example}
  When $n=26$ and $k = 13$, there are 4 equivalence classes of 
  real 2-uniform frames based on conference
  matrices. From Corollary~\ref{cor:556} and
  $1+\sqrt{\frac{(26-13)(26-1)}{13}}=6$, we deduce that the graphs of 
  these frames
  cannot contain any induced complete
  bipartite subgraphs with $m>6$ vertices.
\end{example}

\begin{theorem} Let $F$ be a real 2-uniform $(26,13)$-frame, then
\begin{equation} e_m^{\infty}(F)= \begin{cases}
  \frac{m+4}{10} &\text{if $m \le 6,$} \\
   1        &\text{if $7 \le m \le 26.$} \end{cases} \end{equation}
Consequently, if $F$ is any 2-uniform $(26,13)$-frame, then there
   exists a subset of six frame vectors, $E= \{f_{i_1}, \ldots,
   f_{i_6} \},$
such that $F\backslash E$ no longer spans $\bb R^{13}.$ 
If a set of five or fewer erasures occur, then there exists $L: \bb
   R^{26} \to \bb R^{13}$ such that $LEVx=x$ for all $x \in \bb R^{13}$
   with $\|L\| \le \sqrt{10}.$
\end{theorem}  
\begin{proof} By \cite{BMS} there exist exactly four switching equivalence
   classes of graphs, which by our earlier results give rise to
   exactly four frame equivalence classes of 2-uniform
   $(26,13)$-frames.
Thus, it will sufficient to compute $e_m^{\infty}(F)$ for the frames
   generated by these four graphs. In \cite{BMS} p.~101, representative graphs from
   each of these four equivalence classes are given. A careful
   inspection of these graphs shows that each graph contains a set of
   six vertices such that the induced subgraph is empty and the result
   follows. Combining this fact with our earlier formulas leads to the
   formula for $e_m^{\infty}(F)$ and the estimate on $\|L\|.$
\end{proof}

\begin{remark}
Generally, for a $(26,13)$-frame $F,$ given any set $E$ of 13 or fewer
frame vectors, the set $F\backslash E$ will still span $\bb R^{13}$ and hence
still be a frame. {\rm To see this fact, identify $F$ with its $26 \times
13$ isometric analysis operator $V$. The set of all $26 \times 13$ matrices such
that any collection of 13 or more rows spans $\bb R^{13}$ can easily
seen to be dense in the set of all $26 \times 13$ matrices. If we
polar decompose such a matrix, then it follows that the isometric part
of the polar decomposition inherits this property. Hence, it follows that the
set of $26 \times 13$ isometries such that any set of 13 or more rows
spans $\bb R^{13}$ is dense in
the set of all isometries.}

If we choose any such frame $F_0$ and let $F_1$ denote a 2-uniform frame, 
then $e_m^{\infty}(F_0) <1$ for all $m
\le 13.$ Hence, $e_m^{\infty}(F_0) < e_m^{\infty}(F_1)$ for $6 \le m
\le 13$, while necessarily, $e_2^{\infty}(F_0) \ge e_2^{\infty}(F_1).$
Thus, we see that by minimizing the error for two erasures, we have
necessarily increased the error for some larger number of erasures.
\end{remark}

%\subsection{An exceptional graph}
We continue with an example derived from a graph that is neither of the
conference nor Hadamard type.

\begin{example}
{\it The 2-uniform frame of highest redundancy.}

Among the known graphs giving rise to 2-uniform frames, the frame with
highest redundancy is the 2-uniform $(276,23)$-frame that arises from
the unique regular two-graph on 276 vertices \cite{GS2}. This frame has
redundancy 12. Applying the inequalities of the previous section we
see that
$ e_m^{\infty}(F)
\le \frac{m+4}{60}$ for all $m$ and since the graph cannot contain any
induced complete bipartite subgraphs on 23 or more vertices, this inequality must be
strict for $m>23$. From this formula it follows that if any set of 55
or fewer erasures occurs, then a left inverse for $EV$, $L$, can be constructed
with $\|L\| \le \sqrt{60}.$

>From these inequalitites it follows that given any subset $E$ of $F$
containing at most 56 frame vectors, the set $F\backslash E$ will
still span $\bb R^{23}$. Since we do not precisely know the value of
$e_m^{\infty}(F)$ it is possible that this frame can handle much larger sets
of erasures. By comparison, if we had produced a frame by
simply repeating an orthonormal basis 12 times, then that frame would
be able to handle at most subsets of 11 erasures. On the other hand,
by the argument given in the above remark, a generic uniform
$(276,23)$-frame should be able to handle sets of up to 253 erasures,
but at the expense of having a larger value for $e_2^{\infty}(F).$
\end{example}

We now turn our attention to the special case of graph Hadamards.
Suppose $H$ is a graph Hadamard, which means $H=H^*$, $H^2=nI$, H contains
elements $h_{ij}=\pm 1$ only, and the diagonal is fixed by $h_{jj}=1$. 
The following two results about error 
bounds of frames related to Hadamard matrices are based on an argument 
of Penny Haxell. 

\begin{proposition} \label{prop518}
Any real 2-uniform $(n,k)$-frame $F$ belonging to a signature matrix $Q=H-I$
with a graph Hadamard $H$ satisfies $e_m^\infty(F)=\frac k n +(m-1) c_{n,k}$ 
if $n\ge 48$ and $m\leq 5$.
\end{proposition}
\begin{proof}
By conjugating $H$ by a diagonal matrix of $\pm 1$'s, we may always assume that the first
row and column of $H$ consist entirely of $+1$'s. Then $H^2=nI$ implies that 
the column vectors of $H$ are orthogonal and thus
every additional row and column has to have an equal number of $+1$'s and $-1$'s. 
Moreover, for any two columns other than the first, say $i$ and $j$, their
orthogonality forces them to have $n/4$ entries in common where both are $+1$'s,
$n/4$ entries in common where both columns are $-1$'s, $n/4$ entries
in common where column $i$ is $+1$'s and column $j$ is $-1$'s and
$n/4$ entries in common where column $i$ is $-1$'s and column $j$ is $+1$'s.

The claimed values for the error now follows from showing the existence of
induced complete bipartite subgraphs on $m$ vertices in the
graph $G$ associated with the signature matrix $Q=H-I$. After switching
as described above, we see that $G$ contains one isolated vertex and
that all other vertices have $n/2$ neighbors. So let us pick as $v_1$ the isolated
vertex and as $v_2$ any other, and as $v_3$ one vertex from those $\frac n 2-2$ that are 
neither adjacent with $v_1$ nor $v_2$. Then by the orthogonality argument, $v_3$ has $n/4$ neighbors that 
are adjacent with $v_2$ and a set of $n/4$ neighbors in the set of $\frac n 2 - 2$ that are not adjacent 
with $v_1$ or $v_2$. Thus, there remains a set $A$ of $\frac{n}{4}-3$ vertices that are
not adjacent with $v_1$, $v_2$, or $v_3$. If in this set
there is a pair of vertices that are not adjacent then we have found an induced
edgeless subgraph on 5 vertices, the one consisting of $v_1$, $v_2$,
$v_3$ and the additional two non-adjacent vertices in $A$. Thus, before switching, the subgraph induced by these vertices was complete bipartite.

If there is no non-adjacent pair in $A$, then the induced subgraph on
the vertices in $A$ is a complete graph.
We want to argue that this is impossible for $n$ sufficiently large. Note  that $H$  has eigenvalues $\pm\sqrt n$ because
$H^2=nI$. If $Q$ contains an induced complete subgraph of $s=\frac{n}{4}-3$ vertices,
then the associated signature matrix $Q_s=I_s-J_s$ has eigenvalues $1$ and $1-s$ and 
these must lie between those of $Q$,
$-\sqrt{n}-1\leq 1-s\leq 1\leq\sqrt{n}-1$. Thus, it is impossible that $A$ induces
a complete subgraph  in $G$ when $n>28+8\sqrt{6}>47.5$.
\end{proof}
 
\begin{proposition} \label{prop519}
Any real 2-uniform $(n,k)$-frame $F$ belonging to a signature matrix $Q=I-H$
with a graph Hadamard $H$ satisfies $e_m^\infty(F)=\frac k n +(m-1) c_{n,k}$ 
if $n\ge 30$ and $m\leq 5$.
\end{proposition}
\begin{proof}
The first steps of the proof parallel the one for the preceding proposition,
the only difference being that after switching to obtain the isolated vertex $v_1$,
the valency of the other vertices is $\frac n 2 -2$. Having chosen a vertex $v_2$
and a vertex $v_3$ that is not adjacent with $v_1$ or $v_2$, we observe that
there remains a set $A$ of $s=\frac{n}{4}-2$ vertices that are not adjacent with
any of $v_1$, $v_2$, and $v_3$. As before, we obtain that $Q_s=I_s-J_s$
has eigenvalues $1$ and $1-s$ if $A$ induces a complete subgraph of $s$ vertices in $G$, 
and then necessarily $1-\sqrt{n}\leq 1-s$, thus
there cannot be such an induced complete subgraph if $n>16 +8\sqrt{3}>29.8$.
\end{proof}

\begin{remark}
The smallest possible values of $n$ for graph Hadamards are $n\in\{4, 16, 36, 64\}$.
The preceding results imply that the graphs related to 2-uniform $(n,k)$-frames 
of Hadamard type are guaranteed to contain induced complete bipartites on 5 vertices 
for $n=36,k=15$ and any $n\ge 64$. 
\end{remark}

%List here for page layout!!!!
\begin{example} {\it The 227 known switching equivalence classes of  
graph Ha\-da\-mards with $n=36$.}

For $n=36$, $k=21$, the argument in Proposition~\ref{prop518} does not
guarantee the existence of induced complete bipartite subgraphs on
5 vertices. However, by having a computer search all $227$ known 
equivalence classes \cite{Sp}, one finds that all 
members have at least one induced complete bipartite subgraph on 6 vertices. 
Thus, the $m$-deletion error is the same for all 2-uniform $(36,21)$-frames, 
$$e_m^\infty(F)=\begin{cases}\frac{7+ (m-1)}{12} &\text{ if $m\leq
    6$}\\
1 &\text{if $m > 6$.}\end{cases}$$
Thus, if five or fewer erasures occur, then there exists a left
    inverse of norm at most $\sqrt{12}.$

\begin{table}\label{table2}
\caption{Signature matrices for ``good'' 2-uniform $(36,15)$-frames.}
\begin{minipage}{4.8cm}
\tiny\baselineskip1.5pt\begin{verbatim}
0+++++++++++++++++++++++++++++++++++
+0++++++++++++++++++----------------
++0+++++++++--------++++++++--------
+++0++++++++----------------++++++++
++++0+++----++++----++++----++++----
+++++0--++--++--++--++--++--++--++--
+++++-0-+-+-+-+-+-+-+-+-+-+---++--++
+++++--0-++--+-+-+-+--++--+++-+-+-+-
++++-++-0--+-++-+--+-++-+--+----++++
++++-+-+-0-++--+-++-----++++-++-+--+
++++--++--0+----+++++--+-++-+--+-++-
++++----+++0--++--++-+-+-+-+-+-+-+-+
++--+++--+--0+-++-+++-+--++--++--+-+
++--++-++---+0+--++++-++---+-+--+++-
++--+-+-+--+-+0++++--+-++--+-+++--+-
++--+--+-+-++-+0++-+-++--+-++-++---+
++---++-+-+-+-++0+-+-+--+++-+-+--++-
++---+-+-++--++++0+-+---++-++--++-+-
++----+--++++++--+0++--++-+--+-++--+
++-----++-++++-++-+0-+++--+-+---++-+
+-+-+++---+-++---++-0-++++-++--+-+-+
+-+-++--+--+--+++--+-0+++++-++-++---
+-+-+-+++---++-+---+++0--+++---++-++
+-+-+--+--++-++---++++-0+-+++++--+--
+-+--++-++----+-+++-++-+0-+++-+-+--+
+-+--+---++++--+++--+++--0++-+-+-++-
+-+---++-++-+---+-++-+++++0--++-+-+-
+-+----+++-+-+++-+--+-++++-0--+--+++
+--+++-+--+----+++-+++-++---0-++++-+
+--+++---+-++++---+--+-+-++--0+++++-
+--++-++-+--+-+++------++-++++0--+++
+--++-+---++--++-++-+++--+--++-0+-++
+--+-+-+++---+---+++-++-+-+-++-+0-++
+--+-+--+-++++--+--++--+-+-++++--0++
+--+--+++-+--++-++----+--+++-+++++0-
+--+--+-++-++--+--+++-+-+--++-++++-0
\end{verbatim}
\end{minipage}
\begin{minipage}{4.8cm}
\tiny\baselineskip1.5pt\begin{verbatim}
0+++++++++++++++++++++++++++++++++++
+0++++++++++++++++++----------------
++0+++++++++--------++++++++--------
+++0++++++++----------------++++++++
++++0+++----++++----++++----++++----
+++++0--++--++--++--++--++--++--++--
+++++-0-+-+-+-+-+-+-+-+-+-+---++--++
+++++--0-++--+-+-+-+--++--+++-+-+-+-
++++-++-0--+-++-+--+-++-+--+----++++
++++-+-+-0-++--+-++-----++++-++-+--+
++++--++--0+----+++++--+-++-+--+-++-
++++----+++0--++--++-+-+-+-+-+-+-+-+
++--+++--+--0+-++-+++-+--++--++--+-+
++--++-++---+0+--+++-+++--+--+--+++-
++--+-+-+--+-+0++++--+-++--+-+++--+-
++--+--+-+-++-+0++-+-++--+-++-++---+
++---++-+-+-+-++0+-++---++-++-+--++-
++---+-+-++--++++0+--+--+++-+--++-+-
++----+--++++++--+0++--++-+--+-++--+
++-----++-++++-++-+0+-++---++---++-+
+-+-+++---+-+---+-++0-++++-+++-++---
+-+-++--+--+-+++-+---0+++++-+--+-+-+
+-+-+-+++---++-+---+++0--+++---++-++
+-+-+--+--++-++---++++-0+-+++++--+--
+-+--++-++----+-+++-++-+0-+++-+-+--+
+-+--+---++++--+++--+++--0++-+-+-++-
+-+---++-++-++---++--+++++0---+--+++
+-+----+++-+--+++--++-++++-0-++-+-+-
+--+++-+--+----+++-+++-++---0-++++-+
+--+++---+-++++---+-+--+-+-+-0+++++-
+--++-++-+--+-+++------++-++++0--+++
+--++-+---++--++-++-+++--+--++-0+-++
+--+-+-+++---+---++++-+-+--+++-+0-++
+--+-+--+-++++--+--+-+-+-++-+++--0++
+--+--+++-+--++-++----+--+++-+++++0-
+--+--+-++-++--+--++-++-+-+-+-++++-0
\end{verbatim}
\end{minipage}

\vspace*{0.15cm}

\begin{minipage}{4.8cm}
\tiny\baselineskip1.5pt\begin{verbatim}
0+++++++++++++++++++++++++++++++++++
+0++++++++++++++++++----------------
++0+++++++++--------++++++++--------
+++0++++++++----------------++++++++
++++0+++----++++----++++----++++----
+++++0--++--++--++--++--++--++--++--
+++++-0-+-+-+-+-+-+-+-+-+-+---++--++
+++++--0-++--+-+-+-+--++--+++-+-+-+-
++++-++-0--+-++-+--+-++-+--+----++++
++++-+-+-0-++--+-++-----++++-++-+--+
++++--++--0+----+++++--+-++-+--+-++-
++++----+++0--++--++-+-+-+-+-+-+-+-+
++--+++--+--0+-++-+++-+--++--+-++--+
++--++-++---+0+--+++-++--+-++--++-+-
++--+-+-+--+-+0++++--+-++--++-++---+
++--+--+-+-++-+0++-+-+++--+--++--+-+
++---++-+-+-+-++0+-++--++-+-+---++-+
++---+-+-++--++++0+--+--+++-+-+--++-
++----+--++++++--+0++---++-+-+++--+-
++-----++-++++-++-+0+-++---+-+--+++-
+-+-+++---+-+---+-++0-++++-++++--+--
+-+-++--+--+-+++-+---0+++++--+-+-++-
+-+-+-+++---++-+---+++0--+++--+--+++
+-+-+--+--++--+++--+++-0+-++++-++---
+-+--++-++----+-+++-++-+0-++-++-+-+-
+-+--+---+++++---++-+++--0+++--+-+-+
+-+---++-++-+--+++---+++++0----++-++
+-+----+++-+-++---+++-++++-0+-+-+--+
+--+++-+--+--++-++--+--+-+-+0-++++-+
+--+++---+-++--+--++++-++----0+++++-
+--++-++-+----++-++-+-+-+--+++0--+++
+--++-+---+++++---+--+-+-++-++-0+-++
+--+-+-+++--++--+--+---++-++++-+0-++
+--+-+--+-++---+++-++++--+--+++--0++
+--+--+++-+--+---+++-++-+-+--+++++0-
+--+--+-++-++-+++-----+--++++-++++-0
\end{verbatim}
\end{minipage}
\begin{minipage}{4.8cm}
\tiny\baselineskip1.5pt\begin{verbatim}
0+++++++++++++++++++++++++++++++++++
+0++++++++++++++++++----------------
++0+++++++++--------++++++++--------
+++0++++++++----------------++++++++
++++0+++----++++----++++----++++----
+++++0--++--++--++--++--++--++--++--
+++++-0-+-+-+-+-+-+-+-+-+-+---++--++
+++++--0-++--+-+--++-+-+--+++-+-+-+-
++++-++-0--+-++-+--+-++-+--+----++++
++++-+-+-0-+---++++----++++-+--+-++-
++++--++--0++----++++----+++-++-+--+
++++----+++0--++-+-+--++-+-+-+-+-+-+
++--+++---+-0+-+++-++--++-+--+-++--+
++--++-++---+0+--+++-+++--+-+---++-+
++--+-+-+--+-+0++++--++--+-++-++---+
++--+--+-+-++-+0++-+-+-++--+-+++--+-
++---++-++--+-++0-+++---++-++--++-+-
++---+---+++++++-0+--+--+++--++--+-+
++----++-++--++-++0++-+--++-+-+--++-
++-----++-++++-++-+0+-++---+-+--+++-
+-+-+++---+-+---+-++0-++++-++++--+--
+-+-++-++----+++-+---0+-++++-++-+-+-
+-+-+-+-+--+-++---++++0+-++--+-+-++-
+-+-+--+-+-+++-+---++-+0+-+++--+-+-+
+-+--++-++--+--+++--++-+0-++--+--+++
+-+--+---+++--+-+++-+++--0++++-++---
+-+---++-++-++---++--+++++0----++-++
+-+----++-++--+++--+++-+++-0+-+-+--+
+--+++-+-+---++-+-+-+--+-+-+0-++++-+
+--+++----+++--+-+-++++--+---0+++++-
+--++-++--+---++-++-++--+--+++0--+++
+--++-+--+-++-+++-----++-++-++-0+-++
+--+-+-++-+-++--+--+-+---+++++-+0-++
+--+-+--++-+-+---++++-+++---+++--0++
+--+--++++-----++-++-++-+-+--+++++0-
+--+--+-+-+++++--+-----++-+++-++++-0
\end{verbatim}
\end{minipage}

\vspace*{0.15cm}

\begin{minipage}{4.8cm}
\tiny\baselineskip1.5pt\begin{verbatim}
0+++++++++++++++++++++++++++++++++++
+0++++++++++++++++++----------------
++0+++++++++--------++++++++--------
+++0+++++---+++-----+++-----+++++---
++++0++--++-++-+----+--++---++---+++
+++++0-+-+----++++--++-+-+----++-++-
+++++-0-+--++--+--+++-+--++-+-+--+-+
++++-+-0+-+-+---+++---++-+-++--++-+-
++++--++0+---+--+-++-+-+--+++-+-+--+
+++-++--+0-+-+-+-+-+-+-++-+----++-++
+++-+--+--0++++-+-+---+++-++-+---++-
+++---+--++0-++--++++-+-+++----+++--
++-++-++--+-0-++-+-+--+---++++-+--++
++-++---++++-0+-+-+--+--+---++-+++-+
++-+-+----++++0-++-+++----++-+++-+--
++--+++--+--+--0++++-++-++---++---++
++---+-++-+--+++0-++-+--++-++-+--++-
++---+-+-+-++-++-0+++--+-+-+-+-++--+
++----+++-++-+-+++0---++-+---++-++-+
++----+-++-++-++++-0+---+-+++-+-+-+-
+-+++++----+--+--+-+0--+++-++++-++--
+-++-+--++---++++----0+-++++-+++---+
+-++--++--+++--+--+--+0-+++--++++-+-
+-+-++-++++------++-+--0--++-++-++++
+-+-+----+++-+-++--++++-0+-+++--+-+-
+-+--+++---+---++++-+++-+0-++--+-+-+
+-+---+-+++++-+----+-+++--0+--++-+++
+-+----++-+-+-+-++-+++-++++0++-----+
+--++-+++---++--+--++---++-+0--+++++
+--++-----+-++++-++-+++++--+-0+-+--+
+--+-++-+-----+++-++++++--+--+0-+++-
+--+-+-+-+-++++--+---++--++-+--0++++
+--+---+++-+-+---++++-+++---++++0-+-
+---+++---++-++-+-+-+--+-++-+-++-0++
+---++-+-++-+--++--+--+++-+-+-++++0-
+---+-+-++--++-+-++--+-+-+++++-+-+-0
\end{verbatim}
\end{minipage}
\end{table}
\end{example}

If $n=36$, $k=15$, then Proposition~\ref{prop519} shows that each graph
$G_F$ contains an induced complete bipartite on $5$ vertices. Moreover,
an explicit search finds that
the maximal number of vertices that induce a complete bipartite subgraph
varies from 6 to 8 among the 227 switching-equivalent classes: 
There are 217 switching-equivalent classes that have an induced complete 
bipartite subgraph on 8 vertices, 5 classes that have one on 7 vertices 
but not on 8, and 5 classes that have one on 6 but not on more than 6 vertices.
Thus, for the group of 217, we have that
$e_m^{\infty}(F)=\frac{5+(m-1)}{12}$ for $m \le 8$ and
    $e_m^{\infty}(F) = 1$ for $m \ge 8$.
For the next group of 5 equivalence classes that have an induced 7-vertex complete 
bipartite subgraph, we have 
$e_m^{\infty}(F)=\frac{5+(m-1)}{12}$ if $m \le
    7$, while $e_8^{\infty}(F)\leq \frac{7}{12}+\frac{1}{12}\sqrt{21}\approx 0.965$.
The last bound follows from Theorem~\ref{cor555}.
Finally, for those having a maximal number of 6 vertices that
induce a complete bipartite subgraph, 
we have $e_m^{\infty}(F)= \frac{5+(m-1)}{12}$ for $m \le
    6$,   $e_7^\infty(F)=\frac{13}{24}+\frac{1}{24}\sqrt{65}\approx 0.878$
and $e_8^{\infty}(F) \leq \frac{7}{12}+\frac{1}{12}\sqrt{21}\approx 0.965$.
Again, the results for the cases $m=7$ and $m=8$ follow from
Theorem~\ref{cor555}, because for 
each member of the group, one finds induced subgraphs on 7 vertices that
differ from complete bipartites by only one edge, and we know that
there are no induced complete bipartites on 8 vertices.
The induced subgraphs giving the largest 8-deletion error 
are all found to be switching equivalent and related to complete bipartites by 
flipping two edges. Accordingly, the numerical value for the error 
$e_8^\infty(F)\approx 0.927$ for the members of this group
is below the error bound derived from the absence of complete bipartites.

Thus, if 7 or fewer erasures occur, we know that a left inverse of $EV$
with
norm at most $2\sqrt{6}/\sqrt{11-\sqrt{65}}\approx 2.86$ exists,
compared to $\sqrt{12}\approx 3.46$ for the other 222 switching-equivalent classes.
If 8 erasures occur, we know a left inverse exists of norm at 
most $\sqrt{12}/\sqrt{5-\sqrt{21}}\approx 5.36$.

To summarize, any 2-uniform $(36,15)$-frame 
belonging to the last group of 5 equivalence classes is somewhat preferable 
to the
other 222, because it will have smaller error bounds, but we cannot
guarantee that it can handle any more than 8 erasures. 

We list a representative of the signature matrices belonging to 
each of these ``good'' 5 equivalence  classes in Table~2.

\paragraph{\bf Acknowledgments}
This research was partially supported by the following grants: 
NSF grant DMS-0300128, University of Houston TLCC funds, 
and the Texas Higher Education Coordinating Board Grant 
ARP-003652-0071-2001.

%%%%%%%%%%%%%%%%%%%%%%%%%%%%%%%%%%%%%%%%%%%%%%%%%%%%%%%%%%%%%%%%%%%%%%%%%%%%

%%%

%%%%%%%%%%%%%%%%%%%%%%%%%%%

\end{document}